\documentclass{amsart}

\usepackage{amsmath,amsthm,amsfonts,amscd,amssymb}
\usepackage{url}
\usepackage{comment}

\newtheorem{thm}{Theorem}[section]
\newtheorem{cor}[thm]{Corollary}
\newtheorem{lem}[thm]{Lemma}

\newtheorem{ex}{Example}

\theoremstyle{definition}

\theoremstyle{remark}
\newtheorem{rem}{Remark}[section]

\begin{document}

\title[ Frobenius numbers over totally real number fields]{Positive semigroups and generalized Frobenius numbers over totally real number fields}
\author{Lenny Fukshansky}\thanks{Fukshansky was partially supported by the Simons Foundation grant \#519058}
\author{Yingqi Shi}

\address{Department of Mathematics, 850 Columbia Avenue, Claremont McKenna College, Claremont, CA 91711, USA}
\email{lenny@cmc.edu}
\address{Department of Mathematics, 850 Columbia Avenue, Claremont McKenna College, Claremont, CA 91711, USA}
\email{yshi20@students.claremontmckenna.edu}


\begin{abstract}
Frobenius problem and its many generalizations have been extensively studied in several areas of mathematics. We study semigroups of totally positive algebraic integers in totally real number fields, defining analogues of the Frobenius numbers in this context. We use a geometric framework recently introduced by Aliev, De Loera and Louveaux to produce upper bounds on these Frobenius numbers in terms of a certain height function. We discuss some properties of this function, relating it to absolute Weil height and obtaining a lower bound in the spirit of Lehmer's conjecture for algebraic vectors satisfying some special conditions. We also use a result of Borosh and Treybig to obtain bounds on the size of representations and number of elements of bounded height in such positive semigroups of totally real algebraic integers.
\end{abstract}

\maketitle

\def\A{{\mathcal A}}
\def\AA{{\mathfrak A}}
\def\B{{\mathcal B}}
\def\C{{\mathcal C}}
\def\D{{\mathcal D}}
\def\EE{{\mathfrak E}}
\def\F{{\mathcal F}}
\def\x{{\mathcal H}}
\def\I{{\mathcal I}}
\def\II{{\mathfrak I}}
\def\J{{\mathcal J}}
\def\K{{\mathcal K}}
\def\kk{{\mathfrak K}}
\def\L{{\mathcal L}}
\def\LL{{\mathfrak L}}
\def\M{{\mathcal M}}
\def\mm{{\mathfrak m}}
\def\MM{{\mathfrak M}}
\def\N{{\mathcal N}}
\def\O{{\mathcal O}}
\def\OO{{\mathfrak O}}
\def\PP{{\mathfrak P}}
\def\R{{\mathcal R}}
\def\PNR{{\mathcal P_N(\real)}}
\def\PMNR{{\mathcal P^M_N(\real)}}
\def\PdNR{{\mathcal P^d_N(\real)}}
\def\s{{\mathcal S}}
\def\V{{\mathcal V}}
\def\X{{\mathcal X}}
\def\Y{{\mathcal Y}}
\def\Z{{\mathcal Z}}
\def\H{{\mathcal H}}
\def\BB{{\mathbb B}}
\def\cee{{\mathbb C}}
\def\pee{{\mathbb P}}
\def\que{{\mathbb Q}}
\def\QQ{{\mathbb Q}}
\def\real{{\mathbb R}}
\def\RR{{\mathbb R}}
\def\zed{{\mathbb Z}}
\def\ZZ{{\mathbb Z}}
\def\aaa{{\mathbb A}}
\def\ff{{\mathbb F}}
\def\VV{{\mathbb V}}
\def\kk{{\mathfrak K}}
\def\qbar{{\overline{\mathbb Q}}}
\def\kbar{{\overline{K}}}
\def\ybar{{\overline{Y}}}
\def\kkbar{{\overline{\mathfrak K}}}
\def\ubar{{\overline{U}}}
\def\eps{{\varepsilon}}
\def\ahat{{\hat \alpha}}
\def\bhat{{\hat \beta}}
\def\gt{{\tilde \gamma}}
\def\h{{\tfrac12}}
\def\be{{\boldsymbol e}}
\def\bei{{\boldsymbol e_i}}
\def\bc{{\boldsymbol c}}
\def\bm{{\boldsymbol m}}
\def\bk{{\boldsymbol k}}
\def\bi{{\boldsymbol i}}
\def\bl{{\boldsymbol l}}
\def\bq{{\boldsymbol q}}
\def\bu{{\boldsymbol u}}
\def\bt{{\boldsymbol t}}
\def\bs{{\boldsymbol s}}
\def\bv{{\boldsymbol v}}
\def\bw{{\boldsymbol w}}
\def\bx{{\boldsymbol x}}
\def\bX{{\boldsymbol X}}
\def\bz{{\boldsymbol z}}
\def\bwy{{\boldsymbol y}}
\def\bY{{\boldsymbol Y}}
\def\bL{{\boldsymbol L}}
\def\ba{{\boldsymbol a}}
\def\baa{{\boldsymbol \alpha}}
\def\bb{{\boldsymbol b}}
\def\bet{{\boldsymbol\eta}}
\def\bxi{{\boldsymbol\xi}}
\def\bo{{\boldsymbol 0}}
\def\bone{{\boldsymbol 1}}
\def\bol{{\boldsymbol 1}_L}
\def\ep{\varepsilon}
\def\p{\boldsymbol\varphi}
\def\q{\boldsymbol\psi}
\def\rank{\operatorname{rank}}
\def\aut{\operatorname{Aut}}
\def\lcm{\operatorname{lcm}}
\def\sgn{\operatorname{sgn}}
\def\spn{\operatorname{span}}
\def\md{\operatorname{mod}}
\def\Norm{\operatorname{Norm}}
\def\dim{\operatorname{dim}}
\def\det{\operatorname{det}}
\def\Vol{\operatorname{Vol}}
\def\Area{\operatorname{Area}}
\def\rk{\operatorname{rk}}
\def\ord{\operatorname{ord}}
\def\ker{\operatorname{ker}}
\def\div{\operatorname{div}}
\def\Gal{\operatorname{Gal}}
\def\GL{\operatorname{GL}}
\def\SNR{\operatorname{SNR}}
\def\WR{\operatorname{WR}}
\def\scg{\operatorname{\left< \Gamma \right>}}
\def\swrh{\operatorname{Sim_{WR}(\Lambda_h)}}
\def\ch{\operatorname{C_h}}
\def\cht{\operatorname{C_h(\theta)}}
\def\scgt{\operatorname{\left< \Gamma_{\theta} \right>}}
\def\scgmn{\operatorname{\left< \Gamma_{m,n} \right>}}
\def\gat{\operatorname{\Omega_{\theta}}}
\def\GG{\operatorname{G_{\Lambda_{\ba}}}}
\def\Gg{\operatorname{G}}
\def\Sg{\operatorname{Sg}}
\def\Int{\operatorname{int}}
\def\disc{\operatorname{disc}}

\section{Introduction}
\label{intro}

Let $n \geq 2$ be an integer and let 
\begin{equation}
\label{pos}
1 < a_1 < \dots < a_n
\end{equation}
be relatively prime integers. We say that a positive integer $t$ is {\it representable} by the $n$-tuple $\ba := (a_1,\dots, a_n)$ if
\begin{equation}
\label{rep1}
t = a_1x_1 + \dots + a_nx_n
\end{equation}
for some nonnegative integers $x_1,\dots,x_n$, and we call each such solution $\bx:=(x_1,\dots,x_n)$ of \eqref{rep1} a {\it representation for $t$ in terms of $\ba$}.  Let $s \geq 0$ be an integer, then the {\it $s$-Frobenius number} of this $n$-tuple, $g_s(\ba)$, as defined by Beck and Robins in~\cite{beck_robins}, is the largest positive integer that has at most $s$ distinct representations in terms of $\ba$. This is a generalization of the classical Frobenius number $g_0(\ba)$, i.e., the largest positive integer that has no such representations. The Frobenius number has been studied extensively by a variety of authors, starting as early as late 19th century; see \cite{frobenius_book} for a detailed account and bibliography. The condition 
\begin{equation}
\label{gcd1}
\gcd(a_1,\dots, a_n)=1
\end{equation}
implies that $g_s(\ba)$ exists for every~$s$. The algorithmic {\it Frobenius problem}, known to be NP-hard, is to determine $g_0$ (or more generally $g_s$ for $s \geq 1$) given $n$ and the relatively prime $n$-tuple $a_1,\dots,a_n$ on the input. The hardness of this problem in particular implies that no general closed form formulas for the Frobenius numbers exist, sparking interest in upper and lower bounds.

A geometric approach to the classical Frobenius problem has been pioneered in the influential paper of R. Kannan~\cite{kannan}, leading to a polynomial-time algorithm to find the Frobenius number for each fixed~$n$. Bounds on the classical Frobenius number stemming from further geometry of numbers applications have been obtained in~\cite{me_sinai} and~\cite{iskander}. These ideas have also been extended to the more general $s$-Frobenius problem in~\cite{me_achill} and~\cite{aliev_henk_me}. A higher-dimensional analogue of the Frobenius problem has also been considered in the recent years by several authors, notably in~\cite{aliev_henk}, \cite{aliev_henk_linke}, and~\cite{aliev_deloera}.

This note is inspired by the work of Aliev, De Loera and Louveaux~\cite{aliev_deloera}. We use the geometric setup and results of~\cite{aliev_deloera} (described in Section~\ref{aliev}) to define a natural extension of the Frobenius problem to totally real number fields and to give bounds on the $s$-Frobenius numbers in this context. Let $K$ be a totally real number field of degree $d$ over $\que$ with embeddings $\sigma_1,\dots,\sigma_d : K \to \real$. Let
$$\Sigma := (\sigma_1,\dots,\sigma_d) : K \to \real^d$$
be the Minkowski embedding of $K$. Let $\O_K$ be the ring of integers of $K$, and $\O_K^+$ the additive semigroup of totally positive elements in $\O_K$, i.e.
$$\O_K^+ := \left\{ \alpha \in \O_K : \sigma_i(\alpha) \geq 0\ \forall\ 1 \leq i \leq d \right\}.$$
Let $n > d$ and let $\alpha_1,\dots,\alpha_n \in \O_K^+$ be a collection of elements so that
\begin{equation}
\label{spn_cond}
\O_K = \spn_{\zed} \{ \alpha_1,\dots,\alpha_n \}.
\end{equation}
Such a collection always exists, since there exist bases for $\O_K$ in $\O_K^+$. Indeed, let $\omega_1,\dots,\omega_d$ be a $\zed$-basis for $\O_K$, where $\omega_1 = 1$, and suppose it is not in $\O_K^+$. Let
$$M = \left[ \max_{1 \leq i,j \leq d} |\sigma_i(\omega_j)| \right] + 1,$$
where $[\ ]$ stands for integer part. Then $1, M+\omega_2,\dots,M+\omega_d \in \O_K^+$ and it is still a $\zed$-basis for $\O_K$. Write $\baa = (\alpha_1,\dots,\alpha_n)$. Define the semigroup generated by $\baa$ to be
$$\Sg(\baa) := \left\{ \sum_{i=1}^n \alpha_i x_i : \bx \in \zed^n_{\geq 0} \right\},$$
and the rational cone spanned by $\baa$ to be
$$\C_{\que}(\baa) := \left\{ \sum_{i=1}^n \alpha_i x_i : \bx \in \que^n_{\geq 0} \right\}.$$
Then it is clear that $\Sg(\baa) \subseteq \C_{\que}(\baa) \cap \O_K \subseteq \O_K^+$, and $\Sg(\baa)$ is not necessarily equal to $\C_{\que}(\baa) \cap \O_K$. 

\begin{ex} \label{ex1} Indeed, consider for instance the real quadratic field $K = \que(\sqrt{2})$ and take
$$\alpha_1 = 1,\ \alpha_2 = 4 + \sqrt{2},\ \alpha_3 = 6 + 2\sqrt{2}.$$
One easily checks that these three elements are in $\O_K^+$ and 
$$\O_K = \left\{ a + b \sqrt{2} : a,b \in \zed \right\} = \left\{ (a-4b) \alpha_1 + b \alpha_2 : a,b \in \zed \right\}.$$
Then
$$\Sg(\baa) = \{ (x_1 + 4x_2 + 6x_3) + (x_2 + 2x_3) \sqrt{2} : x_1,x_2,x_3 \in \zed_{\geq 0} \}.$$
On the other hand,
$$3+\sqrt{2} = \frac{1}{2} \alpha_3 \in \C_{\que}(\baa) \cap \O_K,$$
but it is clearly not in $\Sg(\baa)$. 
\end{ex}

\noindent
Further, for each $s \geq 1$ let $\Sg_s(\baa)$ be the set of all points $\beta \in \Sg(\baa)$ for which there are at least $s$ distinct points $\bx \in \zed^n_{\geq 0}$ such that $\sum_{i=1}^n \alpha_i x_i = \beta$, then $\Sg(\baa) = \Sg_1(\baa)$.

Now, let $\beta \in \O_K^+$ be in the interior of the cone $\C_{\que}(\baa)$ and take the ray $t\beta$ as $t \in \zed_{\geq 0}$. Shifting the cone $\C_{\que}(\baa)$ along this ray and intersecting it with $\O_K$ we will eventually land in the semigroup~$\Sg(\baa)$ (this observation will follow from our results). In other words, there exists a positive integer $t$ such that
$$\Int \left( t \beta + \C_{\que}(\baa) \right) \cap \O_K \subseteq \Sg(\baa).$$
More precisely, for each $s \geq 1$ we can define
$$g_s(\baa,\beta) := \min \left\{ t \in \zed_{> 0} : \Int \left( t \beta + \C_{\que}(\baa) \right) \cap \O_K \subseteq \Sg_s(\baa) \right\},$$
and let
$$g_s(\baa) := \max \left\{ g_s(\baa,\beta) : \beta \in \Int(\C_{\que}(\baa)) \cap \O_K^+ \right\}.$$
We refer to $g_s(\baa)$ as the {\it $s$-Frobenius number of $\baa$}. 

We can see that when $K=\que$ this construction reduces to the usual $s$-Frobenius numbers. Indeed, if $K=\que$ then $\O_K=\zed$ and $\O_K^+=\zed_{\geq 0}$, and~\eqref{spn_cond} simply means that $\alpha_1,\dots,\alpha_n$ are positive relatively prime integers. Then $\Sg(\ba)$ is the semigroup of all positive integers representable by $\baa$, $\C_{\que}(\baa) = \que_{\geq 0}$, and so $\C_{\que}(\baa) \cap \zed_{\geq 0} = \zed_{\geq 0}$. Then
\begin{eqnarray*}
& \ & \max_{\beta \in \zed_{\geq 0}} \min \left\{ t \in \zed_{> 0} : \Int \left( t \beta + \que_{\geq 0} \right) \cap \zed \subseteq \Sg_s(\baa) \right\} \\
& \leq & \min \left\{ t \in \zed_{> 0} : \Int \left( t + \que_{\geq 0} \right) \cap \zed \subseteq \Sg_s(\baa) \right\}
\end{eqnarray*}
is precisely the smallest integer $t$ so that all integers $>t$ have at least $s$ representations by~$\baa$.

We present an upper bound on $g_s(\baa)$. To state it, we need to introduce a certain measure of arithmetic complexity of~$\baa$. Let us write $[n] := \{1,\dots,n\}$ and define $\J(n,d) := \{ I \subset [n] : |I| = d \}$. For each $I = \{i_1,\dots,i_d\} \in \J(n,d)$, let us write $\disc(\baa_I)$ for the discriminant of the subcollection $\alpha_{i_1},\dots,\alpha_{i_d}$. We also write $\Delta_K$ for the discriminant of $K$, and define
\begin{equation}
\label{D_measure}
\D(\baa) := \frac{1}{|\Delta_K|} \sum_{I \in \J(n,d)} |\disc(\baa_I)|.
\end{equation}
Notice that absolute values in this definition are not necessary in case of a real number field, since all the quantities are positive; we put them there so that this definition can be naturally extended to any number field. We now state our theorem.

\begin{thm} \label{main1} With notation as above,
$$g_s(\baa) \leq \frac{1}{2\sqrt{n-d+1}} \left( (n-d) \D(\baa) + (s-1)^{\frac{1}{n-d}} \D(\baa)^{\frac{n-d+1}{2(n-d)}} \right).$$
\end{thm}

\noindent
We prove Theorem~\ref{main1} in Section~\ref{aliev}. Now suppose that $\beta \in \Sg(\baa)$, hence there exists $\bx \in \zed^n_{\geq 0}$ such that $\sum_{i=1}^n \alpha_i x_i = \beta$. It is natural to ask for the smallest such representation for~$\beta$. In other words, given $\beta \in \Sg(\baa)$ we want to find $\bx \in \zed^n_{\geq 0}$ such that $\sum_{i=1}^n \alpha_i x_i = \beta$ with $|\bx| := \max \{ |x_i| : 1 \leq i \leq n \}$ as small as possible. This is our next result. To state it, let $\baa(\beta) := (\alpha_1,\dots,\alpha_n, \beta) \in K^{n+1}$ and define
\begin{equation}
\label{M_measure}
\M(\baa,\beta) := \frac{1}{|\Delta_K|^{1/2}} \max_{I \in \J(n+1,d)} |\disc(\baa(\beta)_I)|^{1/2}.
\end{equation}
We also briefly recall the definition of a standard Weil-type height on~$K$. Let us write $M(K)$ for the set of places of~$K$, and for each $v \in M(K)$ let $d_v := [K_v : \que_v]$ be the local degree of $K$ at $v$; in particular, $\sum_{v \mid u} d_v = d$ for each $u \in M(\que)$. Let us normalize absolute values so that the product formula reads
$$\prod_{v \in M(K)} |a|_v = 1$$
for all nonzero $a \in K$. Then the usual {\it inhomogeneous height function} $H_K : K^n \to \real$, which extends Weil height on~$K$ is defined as
\begin{equation}
\label{height_def}
H_K(\baa) = \prod_{v \in M(K)} \max \{ 1, |\alpha_1|_v,\dots,|\alpha_n|_v \}.
\end{equation}
We can now state our next result.

\begin{thm} \label{main2} With notation as above, let $\beta \in \Sg(\baa)$. Then there exists $\bx \in \zed^n_{\geq 0}$ such that $\sum_{i=1}^n \alpha_i x_i = \beta$, and for any such $\bx$ we have 
$$\frac{1}{n} \left( \frac{H_K(\beta)}{H_K(\baa)} \right)^{1/d} \leq |\bx| \leq \M(\baa,\beta).$$
\end{thm}

\noindent
We discuss some properties of $\D(\baa)$ and $\M(\baa,\beta)$ in Section~\ref{heights}, viewing them as kinds of height functions. We use these properties along with Theorem~\ref{main1} to obtain a lower bound on absolute Weil height of~$\baa$ in the spirit of Lehmer's conjecture on heights of algebraic numbers. 

We use a result of Borosh and Treybig~\cite{borosh} to prove Theorem~\ref{main2} in Section~\ref{siegel}. This theorem also allows us to obtain a counting estimate on the number of points of bounded height in the positive semigroup~$\Sg(\baa)$. 

\begin{thm} \label{main3} With notation as above, let $T_1,T_2 \in \real_{\geq 0}$ with $T_1 < T_2$, and define
$$\Sg_s(\baa, T_1, T_2) = \left\{ \beta \in \Sg_s(\baa) : T_1 \leq H_K(\beta) \leq T_2 \right\}.$$
Additionally, for each $\beta \in \Sg(\baa)$ define $r(\beta) = \max \{s : \beta \in \Sg_s(\baa) \}$. Then
$$\sum_{\beta \in \Sg_1(\baa,T_1, T_2)} r(\beta) \leq \left( \frac{d! H_K(\baa)}{|\Delta_K|^{1/2}} T_2 + 1 \right)^n - \left[ \frac{T_1^{1/d}}{n H_K(\baa)^{1/d}} \right]^n.$$
In particular,
\begin{equation}
\label{Sg_s_bnd}
|\Sg_s(\baa,T_1, T_2)| \leq \frac{1}{s} \left\{ \left( \frac{d! H_K(\baa)}{|\Delta_K|^{1/2}} T_2 + 1 \right)^n - \left[ \frac{T_1^{1/d}}{n H_K(\baa)^{1/d}} \right]^n \right\}.
\end{equation}
On the other hand, for $T \geq 1$,
$$\sum_{\beta \in Sg_1(\baa,1, T)} r(\beta) \geq \left[ \frac{T_1^{1/d}}{n H_K(\baa)^{1/d}} \right]^n.$$
\end{thm}

\noindent
Theorem~\ref{main3} is also proved in Section~\ref{siegel}. It is instructive to compare the bounds of Theorem~\ref{main3} to the known estimates on the number of algebraic integers of bounded height in a fixed number field. A result attributed to S. Lang (see~\cite{widmer} for details) asserts that in our case of a totally real number field $K$,
$$\left|\{ \beta \in \O_K : H_K(\beta) \leq T \right\}| = O \left( T^{d^2} (\log T)^{d-1} \right).$$
This implies that our bound~\eqref{Sg_s_bnd} is nontrivial when $n \leq d^2$. We are now ready to proceed.

\bigskip

\section{Polyhedral semigroups and proof of Theorem~\ref{main1}}
\label{aliev}

We start by briefly describing the setup and some results of~\cite{aliev_deloera}. Let $A$ be a $d \times n$ integer matrix, and for each set $I \in \J(n,d)$ let $A_I$ be the $d \times d$ submatrix of $A$ whose columns are indexed by~$I$. Assume that
\begin{enumerate}
\item $\gcd( \det(A_I) : I \in \J(n,d)) = 1$,
\item $\{ \bx \in \real^n_{\geq 0} : A\bx = \bo \} = \{ \bo \}$.
\end{enumerate}
Define the additive semigroup
$$\Sg(A) := \left\{ \bb \in \zed^d : \bb = A\bx \text{ for some } \bx \in \zed^n_{\geq 0} \right\},$$
and for each $s \geq 1$ let $\Sg_s(A)$ be the set of all points $\bb \in \Sg(A)$ for which there are at least $s$ distinct points $\bx \in \zed^n_{\geq 0}$ such that $A\bx = \bb$. Thus $\Sg(A) = \Sg_1(A)$. Let
$$\C_{\real}(A) = \left\{ A \bx : \bx \in \real^n_{\geq 0} \right\}$$
be the convex polyhedral cone spanned by the column vectors of $A$, then it is clear that
$$\Sg(A) \subseteq \C_{\real}(A) \cap \zed^d,$$
and this containment is often proper, i.e. in general $\Sg(A) \neq \C_{\real}(A) \cap \zed^d$. For each $\bb \in \Int(\C_{\real}(A)) \cap \zed^d$, define
$$g_s(A,\bb) = \min \left\{ t \in \zed_{> 0} : \Int \left( t \bb + \C_{\real}(A) \right) \cap \zed^d \subseteq \Sg_s(A) \right\},$$
and let
$$g_s(A) := \max \left\{ g_s(\baa,\beta) : \beta \in \Int(\C_{\real}(A)) \cap \zed^d \right\}.$$
This is the $s$-Frobenius number of $A$ as defined in~\cite{aliev_deloera}. The upper bound on~$g_s(A)$ proved in~\cite{aliev_deloera} (Theorem~2) is
\begin{equation}
\label{A_bnd}
g_s(A) \leq \frac{1}{2\sqrt{n-d+1}} \left( (n-d) \det (AA^{\top}) + (s-1)^{\frac{1}{n-d}}\det (AA^{\top})^{\frac{n-d+1}{2(n-d)}} \right).
\end{equation}
We can now use this result to prove our Theorem~\ref{main1}. Our strategy is straight-forward: we use the Minkowski embedding to convert our setup into that of a lattice in a Euclidean space, prove that our resulting ingredients satisfy the hypotheses of Theorem~2 of~\cite{aliev_deloera}, apply their bound~\eqref{A_bnd}, and then re-interpret it in terms of the original setup in the number field.
\bigskip

\proof[Proof of Theorem~\ref{main1}]
Let the setup be as in Section~\ref{intro}. Let $\alpha_1,\dots,\alpha_n \in \O_K^+$ be a collection of elements satisfying~\eqref{spn_cond}, and let us fix $\omega_1,\dots,\omega_d$, a $\zed$-basis for~$\O_K$. Then there exist integers $a_{ij}$, where $1 \leq i \leq n$, $1 \leq j \leq d$ so that
$$\alpha_i = \sum_{j=1}^d a_{ij} \omega_j.$$
Let us write $A = (a_{ij})^{\top}$ for the $d \times n$ matrix of these integer coefficients. Let
$$B = \begin{pmatrix} \Sigma(\omega_1) & \dots & \Sigma(\omega_d) \end{pmatrix},$$
then $\Delta_K = \det(B)^2$ and
$$C := \begin{pmatrix} \Sigma(\alpha_1) & \dots & \Sigma(\alpha_n) \end{pmatrix} = B A.$$
Since $\alpha_1,\dots,\alpha_n$ satisfy~\eqref{spn_cond}, we must have $\Sigma(\O_K) = C \zed^n$. On the other hand, certainly $\Sigma(\O_K) = B \zed^d$, hence $B \zed^d = B (A \zed^n)$, which means that $A \zed^n = \zed^d$. This implies that row vectors of $A$ are extendable to a basis for $\zed^n$. By Lemma~2 on p.15 of~\cite{cass:geom}, this is equivalent to the condition that
$$\gcd( \det(A_I) : I \in \J(n,d)) = 1.$$
Now suppose $\bx \in \real^n_{\geq 0}$ and assume $A\bx = \bo$. Then
$$C \bx = B (A\bx) = \bo,$$
but entries of $C$ are of the form $\sigma_i(\alpha_j)$, which are all positive real numbers, since $\alpha_j \in \O_K^+$. Therefore $\bx$ must be equal to $\bo$, and so
$$\{ \bx \in \real^n_{\geq 0} : A\bx = \bo \} = \{ \bo \}.$$
Thus matrix $A$ satisfies conditions (1) and (2) above, and so we can apply~\eqref{A_bnd} to get a bound on $g_s(A)$. 

Now notice that $\Sigma(\C_{\que}(\baa)) = B \C_{\que}(A)$, where
$$\C_{\que}(A) := \left\{ A\bx : \bx \in \que_{\geq 0}^n \right\},$$
and $\Int(\C_{\real}(A)) \cap \zed^d = \Int(\C_{\que}(A)) \cap \zed^d$. Indeed, it is clear that 
$$\Int(\C_{\que}(A)) \cap \zed^d \subseteq \Int(\C_{\real}(A)) \cap \zed^d,$$
so let us show containment in the opposite direction. Suppose $\bz \in \Int(\C_{\real}(A)) \cap \zed^d$, then there exists $\bx \in \real^n_{\geq 0}$ such that 
$$A\bx = \bz.$$
In fact, this equation defines a hyperplane in $\real^n$, which is defined over~$\que$ (since $A$ and $\bz$ have integer coordinates), and hence points with rational coordinates are dense in it. Thus taking a sufficiently small open ball in this hyperplane centered at $\bx$, we can find a rational point with positive coordinates satisfying the same equation. This means that $\bz \in \Int(\C_{\que}(A)) \cap \zed^d$.

With this setup in mind, let $\beta \in \Int(\C_{\que}(\baa)) \cap \O_K^+$ and let $\bb \in \zed^d$ be such that $\Sigma(\beta) = B \bb$. Then for $t \in \zed_{> 0}$ we have: 
\begin{eqnarray*}
&\ & \Int \left( t \beta + \C_{\que}(\baa) \right) \cap \O_K \subseteq \Sg_s(\baa) \\
& \Longleftrightarrow & \Int \left( t B\bb + B \C_{\que}(A) \right) \cap B\zed^d \subseteq B \Sg_s(A) \\
& \Longleftrightarrow & \Int \left( t \bb + \C_{\real}(A) \right) \cap \zed^d \subseteq \Sg_s(A).
\end{eqnarray*}
This implies that $g_s(\baa) = g_s(A)$, and so we only need to express~$\det(AA^{\top})$ in terms of~$\baa$. Notice that
$$\det(AA^{\top}) = \det \left( (B^{-1} C) (B^{-1} C)^{\top} \right) = \frac{1}{\det(B)^2} \det(CC^{\top}).$$
Now, $\det(B)^2 = \Delta_K$, and by the Cauchy-Binet formula
$$\det(CC^{\top}) = \sum_{I \in \J(n,d)} \det(C_I)^2,$$
where for each $I = \{ i_1,\dots,i_d \}$,
$$\det(C_I)^2 = \det \begin{pmatrix} \Sigma(\alpha_{i_1}) & \dots & \Sigma(\alpha_{i_d}) \end{pmatrix}^2 = \disc(\baa_I).$$
Combining these observations with~\eqref{A_bnd} completes the proof.
\endproof

\bigskip

\section{Height functions}
\label{heights}

In this section we study some basic properties of the functions~$\D(\baa)$ and $\M(\baa,\beta)$ that we introduced in~\eqref{D_measure} and~\eqref{M_measure}, respectively. Since we view these functions as certain measures of arithmetic complexity, it makes sense to compare them to a traditional height function on number fields. 

Until further notice, let $K$ be any number field of degree $d$, not necessarily totally real as above. As in~\eqref{height_def} above, we write $H_K$ for the inhomogeneous height on~$K$. We can also define the absolute version of Weil height by $H(\baa) = H_K(\baa)^{1/d}$: this height no longer depends on the field of definition. Let us establish some basic properties of~$\D$ and~$\M$ as defined in~\eqref{D_measure} and~\eqref{M_measure}, respectively.

\begin{lem} \label{D_height} Let $K$ be a number field of degree $d$ over $\que$, $n \geq d$, and let $\baa \in K^n$, $\beta \in K$. Then the following are true:
\begin{enumerate}

\item $\D(\baa) = 0$ if and only if $\spn_{\que} \baa \neq K$, and $\M(\baa,\beta) = 0$ if and only if $\spn_{\que} \baa(\beta) \neq K$.

\item If $\baa \in \O_K^n$, then either $\D(\baa) = 0$, or $\D(\baa) \geq 1$; similarly, if $\beta \in \O_K$, then either $\M(\baa,\beta) = 0$, or $\M(\baa,\beta) \geq 1$. Furthermore, $\D(\baa), \M(\baa,\beta) \in \zed_{\geq 0}$.

\item $\D(\baa) \leq \frac{(d!)^2}{|\Delta_K|} \binom{n}{d} H_K(\baa)^2$, $\M(\baa, \beta) \leq \frac{d!}{|\Delta_K|^{1/2}} H_K(\baa) H_K(\beta)$.

\end{enumerate}

\end{lem}

\proof To prove (1), notice that $\D(\baa) = 0$ if and only if discriminant of every $d$-tuple of coordinates of~$\baa$ is equal to~$0$. This happens  if and only if every $d$-tuple of coordinates of~$\baa$ is linearly dependent over~$\que$, meaning that $\spn_{\que} \baa \neq K$. Similarly, $\M(\baa,\beta) = 0$ if and only if discriminant of every $d$-tuple of coordinates of the vector~$\baa(\beta)$ is equal to~$0$, which happens if and only if $\spn_{\que}  \baa(\beta) \neq K$.
\smallskip

To prove (2), assume that $\baa \in \O_K^n$ and $\D(\baa) \neq 0$. Then there exists some $I = \{i_1,\dots,i_d\} \in \J(n,d)$ such that $\disc(\baa_I) \neq 0$. Since coordinates of $\baa$ are algebraic integers, it must be true that for each $1 \leq j \leq d$,
$$\alpha_{i_j} = u_{j1} \omega_1 + \dots + u_{jd} \omega_d,$$
where $\omega_1,\dots,\omega_d$ is a $\zed$-basis for $\O_K$ and $u_{jk}$ are integers such that the matrix $U = (u_{jk})_{1 \leq j,k \leq d}$ is nonsingular. Therefore $\disc(\baa_I) = \Delta_K \det (U)^2$, and hence
$$\D(\baa) \geq \frac{|\disc(\baa_I)|}{|\Delta_K|} =  \det (U)^2 \in \zed_{> 0}.$$
The argument for $\M(\baa,\beta)$ is analogous, replacing $\baa$ with $(\baa,\beta)$ and $\J(n,d)$ with $\J(n+1,d)$, and then observing that $\M(\baa,\beta) \geq \left| \frac{\disc( \baa(\beta)_I)}{\Delta_K} \right|^{1/2} \in \zed_{> 0}$. In particular, $\D(\baa), \M(\baa,\beta) \in \zed_{\geq 0}$.
\smallskip

To prove (3), let $I = \{ i_1,\dots,i_d \} \in \J(n,d)$ and consider the matrix
$$\Sigma(\baa_I) := \begin{pmatrix} \sigma_1(\alpha_{i_1}) & \dots & \sigma_1(\alpha_{i_d}) \\ \vdots & \ddots & \vdots \\ \sigma_d(\alpha_{i_1}) & \dots & \sigma_d(\alpha_{i_d}) \end{pmatrix}.$$
The archimedean absolute values on $K$ are $v_1,\dots,v_d$, given by $|a|_{v_i} = |\sigma_i(a)|^{d_{v_i}}$ for each $a \in K$. Hence we have
\begin{eqnarray}
\label{disc_bnd}
|\disc(\baa_I)|^{1/2} & = & |\det(\Sigma(\baa_I))| = \left| \sum_{\tau \in S_d} \sgn(\tau) \sigma_1(\alpha_{i_{\tau(1)}}) \cdots \sigma_d(\alpha_{i_{\tau(d)}}) \right| \nonumber \\
& \leq & \sum_{\tau \in S_d} |\sigma_1(\alpha_{i_{\tau(1)}})| \cdots |\sigma_d(\alpha_{i_{\tau(d)}})| \nonumber \\
& \leq & d! \prod_{j=1}^d \max \left\{ 1, |\sigma_j(\alpha_{i_1})|^{d_{v_j}},\dots,|\sigma_j(\alpha_{i_d})|^{d_{v_j}} \right\} \nonumber \\
& = & d! \prod_{v \mid \infty} \max \{ 1, |\alpha_{i_1}|_v,\dots,|\alpha_{i_d}|_v \},
\end{eqnarray}
where $\sgn(\tau) = \pm 1$ is the sign of the permutation $\tau$. Therefore
\begin{eqnarray*}
\D(\baa) & = & \frac{1}{|\Delta_K|} \sum_{I \in \J(n,d)} |\disc(\baa_I)| \leq \frac{(d!)^2}{|\Delta_K|} \binom{n}{d} \left( \prod_{v \mid \infty} \max \{ 1, |\alpha_1|_v,\dots,|\alpha_n|_v \} \right)^2 \\
& \leq & \frac{(d!)^2}{|\Delta_K|} \binom{n}{d} \left( \prod_{v \in M(K)} \max \{ 1, |\alpha_1|_v,\dots,|\alpha_n|_v \} \right)^2 = \frac{(d!)^2}{|\Delta_K|} \binom{n}{d} H_K(\baa)^2.
\end{eqnarray*}
This gives the desired bound on $\D(\baa)$ in terms of $H_K(\baa)$. Now, in a manner completely analogous to~\eqref{disc_bnd}, we can obtain
\begin{eqnarray*}
|\disc(\baa(\beta)_I)|^{1/2} & \leq & d! \prod_{v \mid \infty} \max \{ 1, |\alpha_1|_v,\dots,|\alpha_n|_v, |\beta|_v \} \\
& \leq & d! \prod_{v \in M(K)} \left( \max \{ 1, |\alpha_1|_v,\dots,|\alpha_n|_v \} \max \{1, |\beta|_v \} \right) \\
& = & d! H_K(\baa) H_K(\beta),
\end{eqnarray*}
and so $\M(\baa, \beta) \leq \frac{d!}{|\Delta_K|^{1/2}} H_K(\baa) H_K(\beta)$. This completes the proof of the lemma.
\endproof
\medskip

Part (3) of Lemma~\ref{D_height} and its proof show that $\D(\baa)$ and $\M(\baa,\beta)$ measure the arithmetic complexity of $\baa$ and $\baa(\beta)$, respectively, at the archimedean places, allowing for a comparison to the more traditional height function~$H_K(\baa)$. Notice, however, that $\D(\baa)$ and $\M(\baa,\beta)$ are different from traditional heights in the sense that they do not take into account information at the non-archimidean places and do not satisfy Northcott's finiteness property, as does~$H_K(\baa)$: for each $c \in \real_{>0}$, the set $\left\{ \baa \in K^n : H_K(\baa) \leq c \right\}$ is finite, but, say, the set $\left\{ \baa \in K^n : \D(\baa) \leq c \right\}$ is not. Well, clearly $\D(\baa) = 0$ for any $\baa$ with $\spn_{\que} \baa \neq K$ (part (1) of Lemma~\ref{D_height} above), but even sets like
$$\left\{ \baa \in (\O_K^+)^n : \spn_{\zed} \baa = \O_K,  \D(\baa) \leq c \right\}$$
do not have to be finite. 

\begin{ex} \label{ex2} Indeed, let for instance $K=\que(\sqrt{2})$, $t \geq 2$ be a rational integer, and
$$\alpha_1 = 1,\ \alpha_2 = t + \sqrt{2},\ \alpha_3 = 2t + 2\sqrt{2} \in \O_K^+.$$
It is easy to see that $\sqrt{2} = \alpha_2 - t \alpha_1$, hence
$$\O_K = \spn_{\zed} \left\{ \alpha_1,\alpha_2,\alpha_3 \right\}$$
for every $t$. On the other hand, 
$$\D(\baa) = \frac{1}{8} (8 + 32 + 0) = 5,$$
also for every nonzero $t \in \zed$. Hence the set 
$$\{ \baa \in (\O_K^+)^3 : \spn_{\zed} \baa = \O_K,  \D(\baa) \leq 5 \}$$
in this case is infinite.
\end{ex}

\begin{rem} One can think of $\D(\baa)$ as Euclidean norm of the vector of Grassmann coordinates of the matrix $\Sigma(\baa)$, normalized by the discriminant~$\Delta_K$. Transpose of such a matrix can be viewed as a basis matrix for a lattice of rank~$d$ in~$\real^n$, and choosing any other basis for this lattice does not change the value of~$\D(\baa)$. Indeed, rewriting Example~\ref{ex2} in the notation of Section~\ref{aliev}, we have $\Sigma(\baa) = C = BA$, where
$$B = \begin{pmatrix} \Sigma(1) & \Sigma(\sqrt{2}) \end{pmatrix} = \begin{pmatrix} 1 & \sqrt{2} \\ 1 & -\sqrt{2} \end{pmatrix},\ A = \begin{pmatrix} 1 & t & 2t \\ 0 & 1 & 2 \end{pmatrix},$$
and so
$$A^{\top} = \begin{pmatrix} 1 & 0 \\ t & 1 \\ 2t & 2 \end{pmatrix} = \begin{pmatrix} 1 & 0 \\ 0 & 1 \\ 0 & 2 \end{pmatrix} \begin{pmatrix} 1 & 0 \\ t & 1 \end{pmatrix}$$
is a basis matrix of the lattice 
$$\spn_{\zed} \left\{ \begin{pmatrix} 1 \\ 0 \\ 0 \end{pmatrix}, \begin{pmatrix} 0 \\ 1 \\ 2 \end{pmatrix} \right\},$$
where $\begin{pmatrix} 1 & 0 \\ t & 1 \end{pmatrix} \in \GL_2(\zed)$ is just a change of basis matrix. More generally, if $\Sigma(\baa) = BA$ then for any matrix $U \in \GL_d(\zed)$ the matrix $C = B (A^{\top} U)^{\top} = B (U^{\top} A)$ is equal to $\Sigma(\baa')$ for some $\baa' \in \O_K$ with $\D(\baa') = \D(\baa)$, hence there are infinitely many such~$\baa'$. Choosing~$U$ carefully, it is easy to construct such examples with~$\baa' \in \O_K^+$. An analogous observation also applies to the function~$\M(\baa,\beta)$.
\end{rem}

The above observations in particular imply that there cannot exist a general lower bound on $\D(\baa)$ in terms of $H_K(\baa)$: if such a bound existed, we would have
$$\left\{ \baa \in K^n : \D(\baa) \leq c_1 \right\} \subseteq \left\{ \baa \in K^n : H_K(\baa) \leq c_2 \right\}$$
for some appropriate constants $c_1,c_2$. It is therefore interesting to obtain lower bounds on~$\D(\baa)$ besides the trivial one in part (2) of Lemma~\ref{D_height}. Hence, in situations when, say, $g_1(\baa)$ is known (with the setup of Section~\ref{intro}), one can conversely think of Theorem~\ref{main1} as providing a lower bound on $\D(\baa)$:
$$\D(\baa) \geq \frac{2 \sqrt{n-d+1}}{n-d} g_1(\baa).$$
For instance, in a situation when $n=d+1$ and $g_1(\baa) \geq 1$ (a rather common situation, as in Example~\ref{ex1}, for instance), we obtain $\D(\baa) \geq 2\sqrt{2}$ (and hence $\D(\baa) \geq 4$, since it is an integer), which is already better than that of part (2) of Lemma~\ref{D_height}. In fact, $g_1(\baa) \geq 1$ whenever $\Sg(\baa) \neq \C_{\que}(\baa) \cap \O_K$, and hence we have the following immediate corollary. 

\begin{cor} \label{g11} Let the notation be as in Theorem~\ref{main1} and assume that $\Sg(\baa) \neq \C_{\que}(\baa) \cap \O_K$. Then
$$\D(\baa) \geq \max \left\{ 1, \frac{2 \sqrt{n-d+1}}{n-d} \right\}.$$
This lower bound is greater than $1$ as long as $n \leq d+4$.
\end{cor}

\noindent
Combining Corollary~\ref{g11} with part (3) of Lemma~\ref{D_height}, we also obtain the following lower bound on the absolute Weil height of~$\baa$.

\begin{cor} \label{g12} Let the notation be as in Theorem~\ref{main1} and assume that $\Sg(\baa) \neq \C_{\que}(\baa) \cap \O_K$. Then
$$H(\baa) \geq \left( \frac{2 |\Delta_K| \sqrt{n-d+1}\ (n-d-1)! }{d!\ n!} \right)^{1/2d}.$$
\end{cor}

\noindent
The lower bound of Corollary~\ref{g12} is $> 1$ when $|\Delta_K|$ is large in comparison to $n$. For instance, in case $K$ is a quadratic number field and $n=3$, we have
$$H(\baa) \geq \left( \frac{|\Delta_K|}{3\sqrt{2}} \right)^{1/4},$$
while $|\Delta_K|$ can be arbitrarily large (for a general totally real field, Minkowski bound guarantees that $|\Delta_K| \geq \left( d^d/d! \right)^2$).

These observations should be viewed in the light of Lehmer's Problem on lower bound for absolute Weil height of algebraic numbers and the great amount of work done in this direction (see~\cite{mm} for detailed information). Lehmer's conjecture dating back to 1933 states that there exists a constant $c > 1$ such that for every algebraic number~$\alpha$ of degree~$d$, $H(\alpha) > c^{1/d}$. While the conjecture is still open, there is a great number of partial results and generalizations in a variety of special cases. Our lower bound on $H(\baa)$ for the special types of algebraic $n$-tuples~$\baa$ is a small contribution in that general direction.

\bigskip

\section{Proofs of Theorems~\ref{main2} and~\ref{main3}}
\label{siegel}

Let us start from reviewing what we may call a ``positive" version of Siegel's lemma as established by Borosh and Treybig~\cite{borosh}. The name Siegel's lemma often refers to results about the size of solutions of systems of linear equations. The particular version we are interested in is concerned with non-negative solutions to inhomogeneous integer linear systems. Let $A$ be a~$d \times n$ integer matrix such that the equation $A\bx = \bo$ has no nonzero solutions $\bx \in \zed_{\geq 0}^n$. Let $\bb \in \zed^d$ and let $(A\ \bb)$ be the augmented $d \times (n+1)$ matrix. Define 
$$\M(A, \bb) = \max \left\{ |\det \left( (A\ \bb)_I \right)| : I \in \J(n+1,d) \right\}.$$
Theorem~4 of~\cite{borosh} asserts that every $\bx \in \zed_{\geq 0}^n$ such that $A\bx = \bb$ satisfies
\begin{equation}
\label{bbnd}
|\bx| \leq \M(A, \bb).
\end{equation}
We now use this result to prove Theorem~\ref{main2}.

\proof[Proof of Theorem~\ref{main2}]
Let $\beta \in \Sg(\baa)$, then for some $\bx \in \zed_{\geq 0}^n$ we have $\sum_{i=1}^n x_i \alpha_i = \beta$. This means that
$$\Sigma(\beta) = \sum_{i=1}^n x_i \Sigma(\alpha_i) = \Sigma(\baa) \bx.$$
Using the same notation as in Section~\ref{aliev} above, we write $C = \Sigma(\baa) = BA$, and so
$$\beta = C\bx = BA\bx = B\bb,$$
where $\bb \in \zed^m$, hence $A\bx = \bb$. Suppose now that for some $\bwy \in \zed_{\geq 0}^n$ such that $A\bwy = \bo$. Then
$$B A \bwy = C \bwy = \Sigma \left( \sum_{i=1}^n y_i \alpha_i \right) = \bo,$$
which implies that $\sum_{i=1}^n y_i \alpha_i = \bo$. Since $\alpha_1,\dots,\alpha_n \in \O_K^+$, this is only possible if $\bwy = \bo$. Then every $\bx \in \zed_{\geq 0}^n$ such that $A\bx = \bb$ satisfies~\eqref{bbnd}, and for each such~$\bx$ we have $\beta = \sum_{i=1}^n \alpha_i x_i$. Observe also that for each $I \in \J(n+1,d)$,
\begin{eqnarray*}
| \det \left( (A\ \bb)_I \right) | & = & | \det \left( (B^{-1} C\ B^{-1} \Sigma(\beta) )_I \right) | \\
& = & \frac{| \det \left( (\Sigma(\baa)\ \Sigma(\beta))_I \right) |}{| \det(B) |} = \left| \frac{\disc(\baa(\beta)_I)}{\Delta_K} \right|^{1/2},
\end{eqnarray*}
which implies that $\M(A, \bb) = \M(\baa,\beta)$. This completes the proof of the upper bound of the theorem.

To establish the lower bound, let $\baa \in K^n$ and $\bx \in \zed^n$, and let $\beta = \sum_{i=1}^n \alpha_i x_i$. We will prove
\begin{equation}
\label{ht_lemma}
H_K(\beta) \leq n^d |\bx|^d H_K(\baa).
\end{equation}
Indeed, if $v \mid \infty$ in $M(K)$, then
$$|\beta|_v \leq \sum_{i=1}^n |\alpha_i|_v |x_i|_v \leq n |\bx| \max \{ |\alpha_1|_v,\dots,|\alpha_n|_v \}.$$
If $v \nmid \infty$ in $M(K)$, then
$$|\beta|_v \leq \max_{1 \leq i \leq n} |\alpha_i|_v |x_i|_v \leq \max \{ |\alpha_1|_v,\dots,|\alpha_n|_v \},$$
since $\bx \in \zed^n \subseteq \O_K^n$, and so $|x_i|_v \leq 1$ for every $v \nmid \infty$ and $1 \leq i \leq n$. Therefore
\begin{eqnarray*}
H_K(\beta) & = & \prod_{v \in M(K)} \max \{1,|\beta|_v \} \\
& \leq & \prod_{v \mid \infty} \left( n |\bx| \max \{ 1, |\alpha_1|_v,\dots,|\alpha_n|_v \} \right) \times \prod_{v \nmid \infty} \left( \max \{ 1, |\alpha_1|_v,\dots,|\alpha_n|_v \} \right) \\
& \leq & n^d |\bx|^d H_K(\baa).
\end{eqnarray*}
The lower bound of Theorem~\ref{main2} now follows from \eqref{ht_lemma}.
\endproof
\medskip

\proof[Proof of Theorem~\ref{main3}]
Let $1 \leq T_1 < T_2$ be real numbers, then for each $\bx \in \zed_{\geq 0}^n$ such that $\beta = \sum_{i=1}^n \alpha_i x_i \in \Sg_1(\baa, T_1, T_2)$, by Theorem~\ref{main2} and Lemma~\ref{D_height} we have
\begin{eqnarray*}
\frac{1}{n} \left( \frac{T_1}{H_K(\baa)} \right)^{1/d} & \leq & \frac{1}{n} \left( \frac{H_K(\beta)}{H_K(\baa)} \right)^{1/d} \\
& \leq & |\bx| \leq \M(\baa, \beta) \leq \frac{d!}{|\Delta_K|^{1/2}} H_K(\baa) H_K(\beta) \leq \frac{d! H_K(\baa)}{|\Delta_K|^{1/2}} T_2.
\end{eqnarray*}
Now, let
$$\zed_+(\baa,T_1,T_2) = \left\{ \bx \in \zed_{\geq 0}^n : \frac{1}{n} \left( \frac{T_1}{H_K(\baa)} \right)^{1/d} \leq |\bx| \leq \frac{d! H_K(\baa)}{|\Delta_K|^{1/2}} T_2 \right\},$$
and notice that 
\begin{equation}
\label{Z_bnd}
\left| \zed_+(\baa,T_1,T_2) \right| \leq \left( \frac{d! H_K(\baa)}{|\Delta_K|^{1/2}} T_2 + 1 \right)^n - \left[ \frac{T_1^{1/d}}{n H_K(\baa)^{1/d}} \right]^n.
\end{equation}
To each $\beta \in \Sg_1(\baa,T_1, T_2)$ there correspond $r(\beta)$ vectors in $\zed_+(\baa,T_1,T_2)$. Therefore
\begin{equation}
\label{Z_bnd1}
|\Sg_1(\baa,T_1, T_2)| \leq \sum_{\beta \in Sg_1(\baa,T_1, T_2)} r(\beta) \leq \left| \zed_+(\baa,T_1,T_2) \right|.
\end{equation}
Further, for each $\beta \in \Sg_s(\baa,T_1, T_2)$ there are at least $s$ distinct $\bx$ in the set $\zed_+(\baa,T_1,T_2)$, and so
\begin{equation}
\label{Z_bnd2}
|\Sg_s(\baa,T_1, T_2)| \leq \frac{1}{s} \left| \zed_+(\baa,T_1,T_2) \right|.
\end{equation}
On the other hand, for $T >1$ and each $\beta \in \Sg_1(\baa, 1, T)$ there are $r(\beta)$ points in the set
$$\left\{ \bx \in \zed_{\geq 0}^n : |\bx| \leq  \frac{1}{n} \left( \frac{T}{H_K(\baa)} \right)^{1/d} \right\},$$
and so
\begin{equation}
\label{Z_bnd3}
\sum_{\beta \in Sg_1(\baa,1, T)} r(\beta) \geq \left[ \frac{T_1^{1/d}}{n H_K(\baa)^{1/d}} \right]^n.
\end{equation}
The theorem now follows upon combining~\eqref{Z_bnd} with~\eqref{Z_bnd1}, \eqref{Z_bnd2} and~\eqref{Z_bnd3}.
\endproof

\bigskip

\noindent
{\bf Acknowledgement:} We thank the referee for a careful reading of our paper and helpful suggestions that improved the quality of presentation.
\bigskip

\bibliographystyle{plain}  
\bibliography{real_nf_frob}    
\end{document}